\documentclass[12pt]{article}
\usepackage{amsmath,epsfig}

\newcommand{\be}{\begin{equation}}
\newcommand{\ee}{\end{equation}}
\newcommand{\bea}{\begin{eqnarray}}
\newcommand{\eea}{\end{eqnarray}}
\newcommand{\beas}{\begin{eqnarray*}}
\newcommand{\eeas}{\end{eqnarray*}}

\renewcommand{\l}{\left}
\renewcommand{\r}{\right}

\newcommand{\id}{{\bf 1}}
\renewcommand{\vector}{\bar}
\newtheorem{theorem}{Theorem}[section]
\newtheorem{proposition}{Proposition}[section]

\newtheorem{nearremark}{{\it Remark}}[section]
\renewcommand{\Box}{\hbox {$\sqcap$ \kern -1em $\sqcup$}} 
\newenvironment{proof}{\vspace*{2ex}\noindent {\em Proof:} }{$\Box$ \\[2ex]}

\newcounter{stokes}

\newenvironment{acknowledgements}{{\em Acknowledgements:}}{}

\begin{document}

\begin{centering}{\Large Energy in Yang-Mills on a Riemann Surface}		\\[2in]

Dana Fine			\\[1ex]

Mathematics Department		\\
University of Massachusetts	\\
North Dartmouth, MA 02747	\\
dfine@umassd.edu 		\\[2ex]
\end{centering}

\begin{abstract}
Sengupta's lower bound for the Yang-Mills action on smooth
connections on a bundle over a Riemann surface generalize to the space
of connections whose action is finite. In this larger space the
inequality can always be saturated. The Yang-Mills critical sets
correspond to critical sets of the energy action on a space of paths.
This may shed light on Atiyah and Bott's conjecture concerning Morse
theory for ${\cal A / G}$.	\\[3ex]
\end{abstract}

\noindent {\bf Index:} 02.90.+p, 02.40.Vh, 11.15.Tk, 11.10.Kk

\newpage

\renewcommand{\thesection}{\Roman{section}}\vspace{-4ex}
\noindent {\bf Running title:} Energy in Yang-Mills on a Riemann Surface

\section{Introduction}

\renewcommand{\baselinestretch}{2}

One approach~\cite{me1,me2}to quantum Yang-Mills on a Riemann surface
of genus $g$ requires rewriting the Yang-Mills action in
terms of the energy of a $2g$-tuple of paths in the symmetry group
$G$ (This assumes $g \geq 1$. For $g = 0$, the energy is that of a
based loop in $G$.) The energy of such paths appears more recently in
Yang-Mills inequalities Sengupta has developed~\cite{sengupta}.

Sengupta considers the space of smooth
connections, grouped into subspaces by certain requirements on
holonomy. For each subspace, there is a loop in $G$
whose energy bounds from below the Yang-Mills action on that subspace.
For appropriate choices of the requirements on holonomies, this lower
bound can be saturated; Yang-Mills connections are precisely those
which saturate this bound.

Uhlenbeck~\cite{uhlenbeck} has shown that, in two dimensions, the
space of connections whose Yang-Mills action is finite contains
discontinuous connections. Theorem~\ref{inequality} below provides a lower
bound for the Yang-Mills action on this larger space. It is analogous
to Sengupta's, but in this space the bound can always be saturated. One
might then suppose that Yang-Mills connections arise when these
saturating connections are also smooth; this is the import of
Proposition~\ref{relation}.

These relations between the Yang-Mills action and the energy of paths
may help answer a question raised in Atiyah and Bott's seminal
work~\cite{ab} on the topology of the moduli space of Yang-Mills
connections; namely, does the Yang-Mills action, which they show to be
equivariantly perfect, in fact define a Morse stratification?
Theorem~\ref{equiv} describes the correspondence between the critical
sets of the Yang-Mills action and those of the energy on the relevant
space of paths, for which there is reason to believe the analytic
issues are more tractable.

\section{The Geometry of ${\cal A / G}_{m}$}

To describe the required energy requires some background on the structure
of the quotient ${\cal A / G}_{m}$ of the space of connections modulo
gauge transformations. Here ${\cal A}$ refers to connections with
finite total curvature on a given $G$-bundle $P$ over a Riemann
surface $\Sigma$, and ${\cal G}_{m}$ refers to the space of gauge
transformations which are the identity at a specified point $m \in
\Sigma$. What follows is an overview of the essential elements;
details are in the references~\cite{me1, me2}. 

Let $D$, a regular $4g$-gon be a fundamental domain for $\Sigma$,
chosen so that $m$ corresponds to the center of $D$. The edges making
up $\partial D$ represent the generators $\{a_i, b_i\}_{i=1}^{g}$ of
$\pi_1(\Sigma)$, and are identified in pairs, with opposite
orientations, as in Figure 1. 

\begin{figure}
\begin{centering}
\epsfig{file=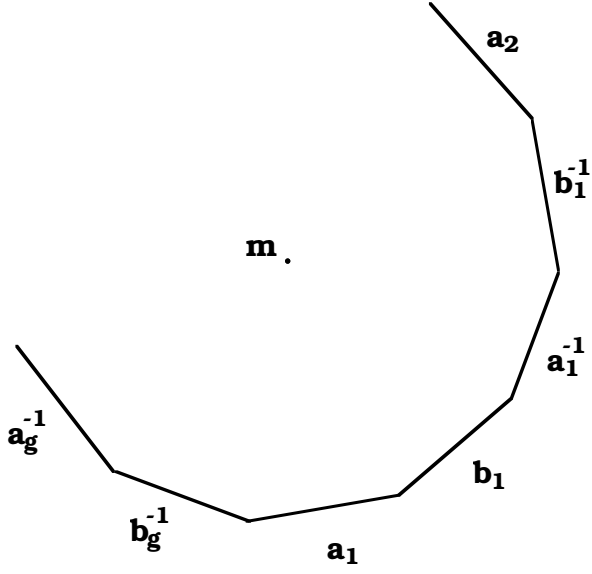}
\caption{The fundamental domain $D$}
\end{centering}
\end{figure}

Theorem 3.1  of~\cite{me2} states that
${\cal A / G}_{m}$ is itself a principal fiber bundle over
$Path^{2g}\,G$ with an affine-linear fiber. Here $Path^{2g}\,G$ is the
space of $2g$-tuples of paths in $G$ subject to a single relation on
the $4g$ endpoint values of the paths. There is an obvious energy
function (see Eq~\ref{energy}) on this base
space $Path^{2g}\,G$; its critical points are are precisely the images of
Yang-Mills connections. To understand how this arises, it will suffice
to examine the projection $\xi: {\cal A / G}_{m} \rightarrow
Path^{2g}\,G$. 

Consider holonomies by a given connection $A$ about the following
loops in $\Sigma$: Pick polar coordinates $(r,\theta)$ on $D$
centered at $m$. For a given point $p$ of the edge $a_1 \subset \partial
D$, the radial path from $m$ to that point followed by the radial path
back to $m$ from the corresponding point $p ^{-1}$ of $a_1 ^{-1}$ defines a
loop in $\Sigma$. See Figure 2. 
\begin{figure}
\begin{centering}
\epsfig{file=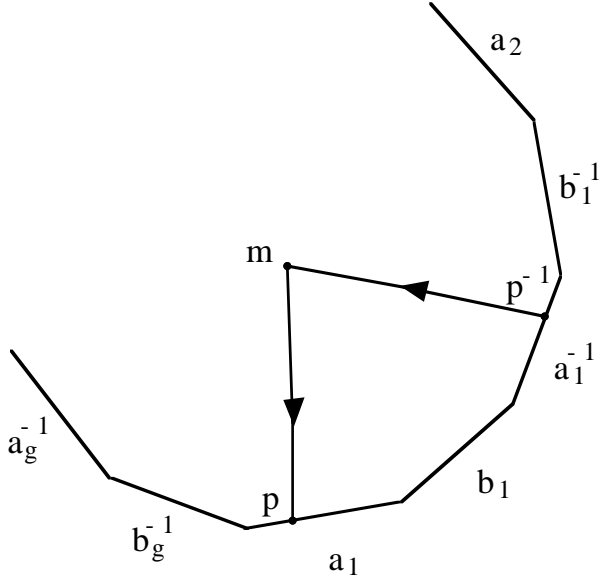}
\caption{A radial path in $D$}
\end{centering}
\end{figure}
Relative to a fixed choice of basepoint in
the fiber over $m$, the holonomy by $A$ about this loop determines an
element of $G$. Now, let the point $p$ vary within $a_1$. The
corresponding holonomies trace out a path $\alpha_1$ in $G$.
Holonomies about radial paths through the points of the other edges
$b_1, a_2, b_2, \cdots, a_g, b_g$ similarly determine paths $\beta_1,
\alpha_2, \beta_2, \cdots, \alpha_g, \beta_g$. Taken together, these
define the $2g$-tuple $\vec{\gamma}_{\! A} = (\alpha_1, \beta_1, \cdots,
\alpha_g, \beta_g)$. These $2g$ paths are not
completely independent of each other, however, as the radii to the
vertices of $\partial D$ each lie on two distinct loops in $\Sigma$
whose holonomies define the endpoint values of distinct paths in $G$.
In fact, traversing, in the appropriate order, each such radius out to
the vertex and back again to $m$ gives a certain product of the
endpoint values of the paths in $\vec{\gamma}_{\! A}$. On the other
hand, by construction, the holonomy about this path must be the
identity in $G$. Equating these gives the relation defining $Path^{2g}\,G$:
\[ 
\alpha_1(0)\beta_1(1)^{-1}\alpha_1(1)^{-1}\beta_1(0) \cdots
\alpha_g(0)\beta_g(1)^{-1}\alpha_g(1)^{-1}\beta_g(0) = \id.
\]

Define $\xi([A]) \equiv \vec{\gamma}_{\! A}$. This is well-defined
on ${\cal A / G}_{m}$, since acting on $A$ by an element of ${\cal
G}_{m}$ has no effect on $\vec{\gamma}_{\! A}$. Clearly, adding to
$A$ a Lie-algebra-valued one-form $\tau$ which vanishes in the radial
directions of $D$ also has no effect on $\vec{\gamma}_{\! A}$. In
fact, in ${\cal A / G}_{m}$, as a bundle over $Path^{2g}\,G$, the
fiber over $\vec{\gamma}_{\! A}$ is the space $\{[A+\tau]: \l. \tau
\r|_{{\scriptstyle \mbox{radii}}} = 0\}$. This, and the fact that $\xi$
is onto, is proven in Theorem~3.1 of~\cite{me2}. 

If the bundle $P$ is not topologically trivial, then, as detailed
in~\cite{me3}, its topology is determined by an element $z$ of the
center of the universal cover $\widehat{G}$ of $G$. On lifting $P$ to a
$\hat{g}$-bundle, the space $Path^{2g}\,G$ is replaced by the
corresponding space for $\widehat{G}$ with the relation
$\alpha_1(0)\beta_1(1)^{-1}\alpha_1(1)^{-1}\beta_1(0) \cdots
\alpha_g(0)\beta_g(1)^{-1}\alpha_g(1)^{-1}\beta_g(0) = z$. Henceforth,
though we omit the hats, we assume we are on the lifted bundle with
the corresponding relation. 

\section{The Yang-Mills action}

Consider now the restriction of the Yang-Mills action on ${\cal A /
G}_{m}$ to the fiber through $\vec{\gamma}_{\! A}$. 
\[
S([A]) = \l< F_{A},  F_{A} \r>,
\]
where the inner product combines the invariant inner product on the
Lie algebra, the metric-induced inner product on forms at each
point and integration over $\Sigma$. Along the fiber, $F_{A + \tau} =
F_{A} + D_{\! A} \tau$, since the term quadratic in $\tau$ vanishes.
Thus,
\[
S([A + \tau]) = S([A]) + 2\l<F_{A}, D_{\! A} \tau\r> + \l< D_{\!
A}\tau,  D_{\! A}\tau \r>.
\]
Theorem~4.2 of~\cite{me2} ensures that the requirement
$\l<F_{\tilde{A}}, D_{\! A}\tau\r>= 0$, for every $\tau$ vanishing
along radii, singles out a unique choice for
a continuous connection $\tilde{A}$ to serve as an ``origin'' in the
fiber. Note that $[\tilde{A}]$ defines a section of ${\cal A / G}_{m}$
over $Path^{2g}\,G$. Relative to this choice of origin,
\begin{equation}
S([\tilde{A}  + \tau]) = S([\tilde{A}]) + \l< D_{\!
A}\tau,  D_{\! A}\tau \r>.		\label{action1}
\end{equation}
(For $\tau$ of the specified form, $D_{\! \tilde{A}} \tau = D_{\!
A}\tau$.)
The key point is that $S([\tilde{A}])$ pulls back to the energy of
$\vec{\gamma}_{\! A}$. This follows from the condition on
$F_{\tilde{A}}$ which implies directly that $*F_{\tilde{A}}$ is
covariantly constant along radii. Thus,
in 
$S([\tilde{A}]) = \l< F_{\tilde{A}},  F_{\tilde{A}} \r>$,
$F_{\tilde{A}}$ may replaced by its average along the radius. This,
however, by a non-Abelian analog of Stoke's theorem, or by Polyakov's
formula, is $\alpha_i ^{-1} \dot{\alpha_i}$ (or $\beta_i ^{-1}
\dot{\beta_i}$) for some $i$ depending on the value of $\theta$. In
fact, for an appropriate choice of parametrization, determined by the
area element on $\Sigma$, 
\begin{equation}
S([\tilde{A}]) = \frac{1}{2} \sum_{i=1}^{2g}\| \dot{\gamma}_i \|^2
\equiv E(\vec{\gamma}_{\! A}),				\label{energy}	
\end{equation}
as detailed in Section~5.1 of~\cite{me2}. Here $\gamma_i$ denotes the
$i$th component of $\vector{\gamma}_{\! A}$. 
For a generic connection, which must be gauge equivalent to $\tilde{A}
+ \tau$, Eq~\ref{action1} thus becomes
\begin{equation}
S([\tilde{A}  + \tau]) = E(\vec{\gamma}_{\! A}) + \l< D_{\!
A}\tau,  D_{\! A}\tau \r>		\label{action}
\end{equation}
It leads immediately to a lower bound
on the Yang-Mills action on a given fiber:
\begin{theorem}
For any connection $A$ representing an element of the fiber through
$\vec{\gamma}_{\! A} \in Path^{2g}\,G$,			\label{inequality}
\[
S([A]) \geq \frac{1}{2}\sum_{i=1}^{2g}\|\gamma_i\|^2,
\]
with equality holding iff $A$ agrees with the section $\tilde{A}$ (up
to gauge transformation).
\end{theorem}
\begin{proof}
This is an immediate consequence of Eq~\ref{action}, since the second
term on the right-hand side is positive semi-definite, and zero iff
$\tau = 0$.
\end{proof}

Given this decomposition of the Yang-Mills action, it is easy to see
how its critical points correspond directly to critical points of the
energy $E$. 
\begin{theorem}		\label{equiv}
The connection $\tilde{A}$ represents a Yang-Mills critical point iff
$\vec{\gamma}_{\! \tilde{A}}$ is a critical point of the energy $E$.
\end{theorem}
\begin{proof}
Suppose $A=\tilde{A} + \tau$ is a Yang-Mills critical point; that is,
a point at which $S([\tilde{A} + \tau])$ is stationary. (There is no
loss of generality in omitting a possible gauge transformation on one
side of this equation.) By considering just $\tau$ of the form $\tau =
t\tau_0$, for $t \in R$, it is clear from Eq~\ref{action} that $\tau = 0$ is a
necessary condition for $A$ to be a critical point. It then follows
that $\vec{\gamma}_{\! A}$ must be a critical point of the energy. The
converse is  immediate.
\end{proof}

To relate this picture, in which connections need not be smooth and
the energy bound can always be saturated, to Sengupta's, in which
connections must be smooth and the energy bound can only be saturated
on the fibers containing Yang-Mills connections, note that in the
fibers over critical points of the energy the connection $\tilde{A}$
is smooth.
\begin{proposition}		\label{relation}
If $\vec{\gamma}$ is a critical point of $E$, then the
corresponding $\tilde{A}$ is smooth.
\end{proposition}
\begin{proof}
A simple calculus of variations computation shows that $\vec{\gamma}$
extremizes $E$ iff 
\[
\frac{\partial \,}{\partial \theta} \gamma_i ^{-1} \dot{\gamma_i} = 0.
\]
On the other hand, this condition also ensures that the covariantly
constant curvatures $F_{\tilde{A}}$, related by the non-Abelian analog
of Stokes Theorem mentioned previously, are continuous at $m$. This
was the only place $\tilde{A}$ might have failed to be smooth.
\end{proof}

\section{A possible application}

Atiyah and Bott suggest equivariant Morse theory might apply to the
cohomology of the Yang-Mills moduli space, and, more particularly,
their stratification may correspond to the Morse stratification for
the Yang-Mills action. With this in mind, they prove the Yang-Mills
action is an equivariantly perfect Morse function. However, analytic
concerns prevent them from developing the theory more fully, except in
genus zero. There Bott and Samuelson~\cite{bs} have shown that ${\cal
A / G}$ is topologically equivalent to based loops in $G$, and that
Morse theory arguments go through for a wide variety of symmetric
spaces including these based loops.

The geometric picture of ${\cal A / G}_{m}$ as an affine-linear bundle
shows it is topologically equivalent to its base space
$Path^{2g}\,G$. Passing from ${\cal G}_{m}$ to ${\cal G}$, this
becomes $Path^{2g}\,G/G$, where a given element $g \in G$ acts
adjointly on each path: $\gamma_i(t) \mapsto g ^{-1} \gamma_i(t) g$.
Moreover Eq~\ref{action} says the section $[\tilde{A}]$ pulls the
Yang-Mills action back to the energy on $Path^{2g}\,G/G$. Clearly, the
Morse theory for this base space, if such exists, would be the Morse
theory for ${\cal A / G}$. Moreover, the generality of Bott and
Samuelson's results is nearly sufficient to apply them directly to
$Path^{2g}\,G/G$. The endpoint condition, however, requires careful
treatment, which we defer to future work.

\begin{acknowledgements} The author is grateful to Ambar Sengupta for
discussions of his work on the energy inequality and to Stephen Sawin
for discussions of this and many other aspects of two-dimensional
Yang-Mills.
\end{acknowledgements}

\end{document}